\documentclass[a4paper,12pt]{article}
\usepackage{amsfonts}
\usepackage{amsthm}
\usepackage{amsmath}
\usepackage{amsfonts}
\usepackage{latexsym}
\usepackage{amssymb}
\usepackage[latin1]{inputenc}

\newtheorem{fed}{Definition}[section]
\newtheorem{teo}[fed]{Theorem}
\newtheorem*{teo*}{Theorem}
\newtheorem{lem}[fed]{Lemma}
\newtheorem{cor}[fed]{Corollary}
\newtheorem{pro}[fed]{Proposition}  
\theoremstyle{definition}
\newtheorem{rem}[fed]{Remark}
\newtheorem{rems}[fed]{Remarks}
\newtheorem{exa}[fed]{Example}

\def\barr{\begin{array}}
\def\earr{\end{array}}

\newcommand{\NN}{\mathbb{N}}

\newcommand{\F}{\mathbb{F}}
\def\NN{\mathbb{N}}
\def\R{\mathbb{R}}
\def\C{\mathbb{C}}

\def\noi{\noindent}
\def\bma{\left[\begin{array}}
\def\ema{\end{array}\right]}
\def\ben{\begin{enumerate}}
\def\een{\end{enumerate}}
\newcommand{\IN}[1]{\mathbb {I} _{#1}}
\newcommand{\ind}[1]{\mbox{Ind} (#1)}
\newcommand{\exe}[1]{E \left(#1\right)}

\newcommand{\frami}{\{f_i\}_{i\in I}}

\usepackage[dvips]{color}

\def\ds{\displaystyle}
\def\noi{\noindent}
\definecolor{color}{RGB}{140,140,140}
\def\QED{\hfill $\color{color}{\blacksquare}$}
\def\EOE{\hfill $\color{color}{\blacktriangle}$}
\def\bdem{\begin{proof}}
\def\bdem{\begin{proof}}

\def\N{\mathcal{N}}

\def\cB{\mathcal{B}}

\def\cE{\mathcal{E}}
\def\cF{\mathcal{F}}
\def\cG{\mathcal{G}}
\def\cK{\mathcal{K}}

\def\cM{\mathcal{M}}
\def\cN{\mathcal{N}}
\def\cP{\mathcal{P}}

\def\cS{\mathcal{S}}

\def\cU{\mathcal{U}}
\def\cW{\mathcal{W}}

\def\pesos{\cP\left(\cW \right)}

\DeclareMathOperator{\sspan}{span}

\newcommand{\pint}[2]{\displaystyle \left \langle #1 \,, \, #2  \right\rangle}

\newcommand{\hil}{\mathcal{H}}
\newcommand{\op}{L(\mathcal{H})}

\newcommand{\posop}{L(\mathcal{H})^+}

\newcommand{\convk}{\xrightarrow[k\rightarrow\infty]{}}

\newcommand{\tfs}{T_{\cW_w}}
\newcommand{\frs}{S_{\cW_w}}
\newcommand{\ra}{\rightarrow}
\def\ssec{\cW = \{ W_i\}_{i\in I}}
\def\sseck{\cW = \{ W_k\}_{k\in \NN}}

\def\refi{\cV = \{ V_i\}_{i\in J}}

\def\FS {\cW_w\,}
\def\sfram{\cW_w  = (w_i\, ,\, W_i)_{i\in I}}

\def\sframN{\cW_w = (w_k\, ,\, W_k)_{k\in \NN}}

\def\inv{^{-1}}

\def\beq{\begin{equation}}
\def\eeq{\end{equation}}
\def\N{\mathbb{N}}
\def\F{\mathcal{F}}

\def\cB{\mathcal{B}}

\def\cE{\mathcal{E}}
\def\cH{\mathcal{H}}
\def\cK{\mathcal{K}}
\def\cU{\mathcal{U}}

\def\cP{\mathcal{P}}

\def\ese{\mathcal{M}}

\def\eme{\mathcal{M}}
\def\ene{\mathcal{N}}
\def\cV{\mathcal{V}}
\def\cW{\mathcal{W}}

\def\inc{\subseteq}

\def\onb{orthonormal basis }

\def\vacio{\varnothing}
\def\orto{^\perp}
\def\inc{\subseteq}
\def\sii{ if and only if }
\def\inv{^{-1}}
\def\*A{\#\sb A}

\def\H{{\cal H}}

\def\glh{Gl(\cH)}

\def\glkw{Gl(\cK_\cW)}
\def\ca{L(\H ) }
\def\cam{L(\H )^+ }

\def\cH{{\cal H}}

\def\cP{{\cal P}}

\def\cS{{\cal S}}

\def\cM{{\cal M}}
\def\cN{{\cal N}}

\def\rai{^{1/2}}
\def\mrai{^{-1/2}}

\def\api{\langle}
\def\cpi{\rangle}

\def\noi{\noindent}

\def\bm{\left(\begin{array}}
\def\em{\end{array}\right)}
\def\ben{\begin{enumerate}}
\def\een{\end{enumerate}}
\def\barr{\begin{array}}
\def\earr{\end{array}}
\def\iiff{if and only if }
\def\inv{^{-1}}

\def\H{{\cal H}}

\def \linm{\ell^\infty _+ (I)}
\def \linmJ{\ell^\infty _+ (J)}
\def \linmN{\ell^\infty _+ (\NN)}

\newcommand{\fram}{\{f_n\}_{n\in \mathbb{N}}}

\newcommand{\bon}{\{e_k\}_{k\in \mathbb{N}}}

\def\lh{{L(\H)}}
\def\lk{{L(\cK)}}
\def\lh+{{\lh^+}}
\def\lhk{L(\H , \cK)}
\def\lkh{L( \cK , \H )}

\def\noi{\noindent}

\newcommand{\dist}[2]{\mbox{dist}\left(#1,\, #2\right)}

\newcommand{\angf}[2]{c\left[\,#1,\,#2\,\right]}
\newcommand{\sang}[2]{s\left[\,#1,\,#2\,\right]}

\def\inc{\subseteq}
\def\bm{\left(\begin{array}}
\def\em{\end{array}\right)}
\def\subcer{\sqsubseteq}
\newcommand{\peso}[1]{ \quad \text{ #1 } \quad }

\newcommand{\trivial}{\{0\}}
\newcommand{\gen}[1]{\overline{\mbox{span}}\left\{#1\right\}}
\def\bsi{ \linm\, ^* }

\date{}
\begin{document}
\title{Frames of subspaces  and operators}
\author{Mariano A. Ruiz$^*$  and Demetrio Stojanoff\footnote{Partially supported by CONICET (PIP 4463/96), Universidad de La Plata (UNLP 11 X350) and
ANPCYT (PICT03-09521).}
}

\maketitle

\begin{abstract}
We study the relationship between operators, orthonormal basis of subspaces and frames of subspaces (also called fusion frames) for a separable Hilbert space $\mathcal{H}$. We get sufficient conditions on an orthonormal basis of subspaces $\mathcal{E} = \{E_i \}_{i\in I}$ of a Hilbert space $\mathcal{K}$ and  a surjective $T\in L(\mathcal{K} , \mathcal{H} )$ in order  that $\{T(E_i)\}_{i\in I}$ is a frame of subspaces with respect to a computable sequence of weights.
We also obtain generalizations of  results in [J. A. Antezana, G. Corach, M. Ruiz and D. Stojanoff,  Oblique projections and frames. Proc. Amer. Math. Soc.  134  (2006), 1031-1037],  which related frames of subspaces (including the computation of their weights) and oblique projections. The notion of refinament of a fusion frame is defined and used to obtain results about the excess of such frames. We study the set of admissible weights for a generating sequence of subspaces. Several examples are given.
\end{abstract}

\vglue1truecm

\noi
{\bf Keywords:} frames, frames of subspaces, fusion frames, Hilbert space operators, oblique projections.

\medskip
\noi
{\bf 2000 AMS Subject Classifications:} Primary 42C15, 47A05.

\section{Introduction}
Let $\hil$ be a (separable) Hilbert space. A {\sl frame} for $\hil$ is a sequence of vectors 
$\F = \frami $ for which  there exist numbers $ A,B>0$ such that
$$
A\|f\|^2\leq \sum_{i \in I}
|\pint{f }{ f_i}|^2 \leq B\|f\|^2 \ , \peso{  for every} f \in \H \ .
$$
This definition has been  generalized to the notion of frames of subspaces by Casazza and Kutyniok
 \cite{[CasKu]} (see also \cite{[For1]} and \cite{[For2]}) in the following way:  
Let  $\ssec$ be a sequence of closed subspaces, and let $w = \{w_i\}_{i \in I} \in \linm $ 
(i.e. $w_i>0$ for every $i\in I$). We say that $\sfram$  
is a frame of subspaces  (shortly: FS) for $\hil$  if there exist $A _{\cW_w}\, ,\ 
B_{\cW_w} > 0$  such that 
\[
A_{\cW_w} \, \|f\|^2\leq \sum_{i \in I}
w_i^2\| P_{W_i}f\|^2 \leq B_{\cW_w} \, \|f\|^2   \peso{for every $f\in \hil$  ,}
\]
where each $P_{W_i}$ denotes the orthogonal projection onto $W_i\,$. 
The relevance of this notion, as remarked in \cite {[CasKu]},  
is that it gives criteria for constructing a frame for $\hil$, by joining 
sequences of frames for subspaces of $\hil$ (see Theorem \ref{CK} for details). 
In other words, to give  conditions  which assure that a sequence of 
``local" frames, can be pieced together to obtain a frame for the complete space. 

Recently, the frames of subspaces have been renamed as \it fusion frames.  \rm 
This notion is 
intensely studied during the last years, and several new applications 
have been discovered. The reader is referred to Casazza, Kutyniok, Li 
\cite {[CasKuLi]}, Casazza and Kutyniok \cite{[CasKu2]}, Gavruta \cite{[Gav]} and the 
references therein.

Given  sequences $\sfram$, consider, for each $i \in I$, 
an \onb $\{e_{ik}\}_{k\in K_i}$  of $W_i\,$. It was proved in \cite{[CasKu]} that $\cW_w$ is a FS 
for $\hil$ \iiff $\cE = \{w_i e_{ik}\}_{i\in I, k\in K_i}$ is a frame for $\hil$. 
Therefore,  a FS can be thought as a frame (of vectors) such that some subsequences 
are required to be orthogonal and to have the same norm. 
Therefore, many objects associated to vector frames have a generalization 
for frames of subspaces (see \cite {[CasKu]} and \cite{[AsKho]}), for example, 
synthesis, analysis and frame operators. 
Also,  some useful results concerning frames  still hold in the FS setting.
For instance, as it is shown in  \cite{[AsKho]}, a Parseval FS is an 
orthogonal projection of a orthonormal basis of subspaces of a larger Hilbert 
space containing $\hil$, generalizing the well known result of D. Han and D. Larson.

As we mention before, if $\sfram$ is a FS, 
the synthesis, 
analysis and frame operator can be defined, and  the properties of $\cW_w$ can be study using these operators, 
as well as for frames of vectors  (\cite {[CasKu]}, \cite{[AsKho]}). In \cite{[CasKu]}, 
the domain of the synthesis operator is defined as $\cK_\cW =\bigoplus_{i\in I}W_i\, $. 
So the subspaces $\{W_i\}_{i\in \N}$ are embedded in $\cK_\cW $ as an orthonormal basis
of subspaces (see also \cite{[AsKho]} where other type of domain is used).  
Therefore, the frame of subspaces is the image of the orthonormal basis under the 
synthesis operator (which is a bounded surjective operator). 
 
 However, if $\sfram$ is a FS for $\H$,  its synthesis operator $T_{\cW_w} $ satisfies that  
$T_{\cW_w} g = w_i \, g$ for every $g $ in the copy of each $W_i$ into  $\cK_\cW$ 
(see \cite{[CasKu]} or Definition \ref{defis} below).  
Hence, unlike   the vector case,  if one fix 
an \onb of subspaces $\cE = \{E_i \}_{i\in I}$ of 
a Hilbert space $\cK$, not every  surjective operator $T\in L(\cK , \H )$ 
is the synthesis operator of a FS. Even worse,   
there exist surjective operators $T\in L(\cK , \H )$   
such that  $T(E_i)$ is closed for every $i\in I$,  but the sequence 
$ (w_i \, , \, T(E_i)\,)_{i \in I}$ fails to be a FS for every 
$w\in \linm \,$ (see Example  \ref{gama no frame}).

The purpose of this work is to study the relationship between operators and frames of subspaces. 
Our aim is to get 
more flexibility in the use of operator theory 
techniques, with respect to the rigid definition of the synthesis operator. In this direction 
we get (sufficient) conditions on an \onb of subspaces $\cE = \{E_i \}_{i\in I}$ of 
a Hilbert space $\cK$  and a surjective $T\in L(\cK , \H )$ in order to assure that 
they produce a frame of subspaces with respect to a computable sequence of weights 
(Theorem \ref{cotas}). We use then this result for describing properties of 
equivalent frames of subspaces, and 
for studying the $excess$ of such  frames. We 
obtain generalizations of two results of \cite{[ACRS2]}, which relate FS 
(including the computation of their weights) and oblique projections
(see also \cite{[AsKho]} and \cite{[CasKuLi]}). 
We also define the notion of refinement of sequences of subspaces and frames 
of subspaces. This allows us to describe the excess of frames of subspaces, obtaining 
results which are very similar  to the known results in classical frame theory. 

It is remarkable that several known results of frame theory are not valid  in the FS setting. 
For example,  we exhibit  a frame of subspaces  $\sfram$ of $\hil$ such that, for every $G \in \glh$, 
the sequence $(v_i  \, , \,  G W_i )_{i \in I}$ fails to be a Parseval FS for every 
$v \in \linm$, including the case  $G = S_{\cW_w}\mrai \,$, where $S_{\cW_w}$ is the frame 
operator of $\FS$ (see Examples \ref{frame no gama} and \ref{otros2}). Several of this facts 
are exposed in a section of (counter)examples.

Finally  we begin with the study of that is, in our opinion,  the key problem of the theory 
of frames of subspaces: given a generating sequence $\ssec$ of closed subspaces of $\hil$, 
to obtain a characterization of the set of its admissible weights, 
\[
\pesos= \big\{ \ w \in \linm : \  \sfram \text { is a FS for } \hil \big\} 
  \ .
\]
Particularly, we search for conditions  which assure that a sequence $\cW$ satisfy that 
$\pesos \neq \vacio$. We obtain some partial results about these problems, and we study 
an equivalent relation between weights, compatible with their admissibility 
with respect to a generating sequence. We give also several examples 
which illustrate the complexity of the problem.

The paper is organized as follows: Section 2 contains preliminary results 
about angles between closed subspaces, the reduced minimum modulus of operators, 
and frames of vectors. In section 3 we introduce the frames of subspaces and 
we state the first results relating these frames and Hilbert space operators. 
In Section 4 the set of admissible weights of a FS is studied. 
Section 5 contains the results which relate oblique projections and frames of subspaces.
Section 6 is devoted to refinement of sequences of subspaces and it contains several
results about the excess of a FS. In section 7 we present a large collection  
of examples.

\medskip

\noindent {\it Note}: after completing this paper, the authors were pointed out of the existence of recent works on fusion frames \cite{[CasKuLi]}, \cite{[CasKu2]} and \cite{[Gav]} . Thus, Corollary \ref{conG}  appears in \cite{[CasKuLi]} and \cite{[Gav]}. Also, Theorem \ref{T:poyeccionoblicua a frames} is related  with Theorem 3.1. in \cite{[CasKuLi]}. Nevertheless, the proofs in general are quite different.

\section{Preliminaries and Notations.}

Let $\hil$ and $\cK$ be separable Hilbert spaces and 
$\lhk $  the space of bounded linear operators $A:\hil \to \cK$ 
(if $\cK = \H$ we write $\ca$\,). The symbol $\glh$ 
denotes the group of invertible operators in $\ca$, and
$\glh^+$ the set of positive definite invertible operators on $\H$. 
For an operator $A \in \lhk$,  $R(A)$ denotes the range of $A$, $ N( A)$ the
nullspace of $A$,  $A^* \in \lkh$ the adjoint
of $A$, and $\|A\|$ the
operator norm of $A$. 

We write $\ese\subcer \hil$ to denote that $\ese$ is a {\bf closed subspace} of $\H$. 
Given $ \ese \subcer \hil$, $P_\ese$ is the
orthogonal (i.e., selfadjoint) projection onto $\ese$. 
If also $\ene\subcer\hil$, we write $\ese \ominus \ene := \ese \cap (\ese \cap \ene) \orto$. 

Let $I$ be a denumerable set. 
We denote by $\linm$ the space of bounded sequences of positive numbers. 
We consider on $\linm$ the usual product of $\ell^\infty  (I)$ (i.e. cordinatewise product). With this product
$\ell^\infty(I) $ is a von Neumann algebra. We denote by 
\beq\label{bsi}
\bsi = \{ \, \{w_i\}_{i\in I} \in \linm : \inf _{i \in I} w_i > 0 \} = \linm \cap Gl (\ell^\infty(I)\,) \ .
\end{equation} 
We shall recall the definition and basic properties of angles between closed subspaces of
$\hil$. We refer the reader to \cite{[ACRS]} for details and proofs. See also the survey by  Deutsch
\cite{[De]} or the book by Kato \cite{[Kato]}.

\subsection*{Angle between subspaces and reduced minimum modulus.}
We shall recall the definition of angle between closed subspaces of
$\hil$. We refer the reader to \cite{[ACRS]} (where the same notations are used)  
for details and proofs. See also the survey by  Deutsch
\cite{[De]} or the book by Kato \cite{[Kato]}.

\begin{fed}\rm
Let $ \eme , \ene \subcer \H$. The \textbf{angle} between $\eme$ and
$\ene$ is the angle in $[0,\pi/2]$ whose cosine is 
\[
\angf{\eme}{\ene}=\sup\{\,|\pint{x}{y}  |:\;x\in \ese \ominus \ene  , \ y\in \ene \ominus \eme \ \mbox{and}\   \|x\|=\|y\|=1 \} \ .
\]
If $\ese \inc \ene$ or $\ene \inc \ese$, we define $\angf{\eme}{\ene} =0$, as if they 
where orthogonal. 
The $sine$ of this angle is denoted by $\sang{\eme }{\ene } = (1-\angf{\eme }{\ene }^2\,)\rai$. \EOE
\end{fed}

\medskip
\noi Now, we state some known results concerning angles  (see  \cite {[ACRS]} or \cite{[De]}).

\begin{pro}\label{propiedades elementales de los angulos} \rm
Let $ \eme , \ene \subcer \H$. Then 
\ben
\item $\angf{\eme}{\ene}=\angf{\ene}{\eme}=\angf{\ese \ominus \ene}{\ene}=
  \angf{\eme}{\ene \ominus \eme}$.
\item If $\dim \ese < \infty$, then $\angf{\eme}{\ene}<1$.
\item $\angf{\eme}{\ene}<1$   \sii $\eme+\ene$ is closed.
\item $\angf{\eme}{\ene}= \angf{\eme^\bot}{\ene^\bot}$
\item $\angf{\eme}{\ene}=\|P_\eme P_{\ene \ominus \eme}\|=\|P_{\ese \ominus \ene}P_\ene \|
=\|P_\eme P_\ene-P_{\eme\cap \ene}\|$.
\item $\sang{\eme}{\ene}= \dist{B_1(\eme \ominus \ene) }{\ene}$, where  $B_1(\eme \ominus \ene)$ 
is the unit ball of $\eme \ominus \ene$.
\QED
\een
\end{pro}

\begin{fed}\rm
The \textit{reduced minimum modulus} $\gamma(T)$ of 
$T\in \lhk$ is defined by
\begin{equation}\label{gamma}
\gamma (T)= \inf \{ \|Tx\| :  \|x\|=1 \; , \;  x\in N(T)^\bot \}
\end{equation}
\end{fed}

\begin{rem}\label{Pgamma}
The following properties are well known (see \cite{[ACRS]}). Let $T \in \lhk$. 
\ben
\item $\gamma (T)=\gamma (T^*) =\gamma (T^*T)\rai $. 
\item $R(T) \subcer \cK$ \sii $\gamma (T)>0$.    
\item If $T$ is invertible, then $\gamma (T) = \|T\inv\|\inv$. 
\item If $B\in \lk$, then 
\begin{equation}\label{desi1}
\|B^{-1}\|^{-1} \gamma (T) \leq \gamma (BT) \leq \|B\| \gamma (T) \ .
\end{equation}
\item Suppose that  $R(T)\subcer \cK $ and take 
$\eme \subcer \hil$. Then
\begin{equation}\label{desi2}
 \gamma (T)\ \sang{N(T)}{ \eme} \leq \gamma (TP_\eme)\leq \|T\| \ \sang{N(T)}{ \eme} .
\end{equation}
 In particular, $T(\ese) \subcer   \cK$ \sii $ \angf{N(T)}{ \eme}<1$. 
\EOE
\een
\end{rem}

\subsection*{Preliminaries on frames.}

We introduce some basic facts about frames in Hilbert spaces. For
a complete description of frame theory and its applications, the reader
is referred to
Daubechies, Grossmann and Meyer \cite {[DGM]}, 
the review by Heil and Walnut \cite{[HeiWal]}
or the books by Young \cite{[Y]} and Christensen
\cite{[liChr]}.
\begin{fed} \rm
Let $\F = \fram$ a sequence in a Hilbert space $\H$.
$\F$ is called a \textsl{frame} if
there exist numbers $ A,B>0$ such that
\beq\label{frame}
A\|f\|^2\leq \sum_{n\in \mathbb{N}}
|\pint{f }{ f_n}|^2 \leq B\|f\|^2 \ , \peso{  for every} f \in \H \ .
\end{equation}
The optimal constants  $A_\F\, , B_\F$ for Eq. (\ref{frame})  are called
the \textsl{frame bounds} for $\F$.
The frame $\F$ is called $tight$ if $A_\F=B_\F$, and \it Parseval \rm  if
$A_\F=B_\F=1$. 
\EOE
\end{fed}

\begin{fed}\rm\label{preframe}
Let $\F= \fram$ be a frame in $\hil$ and let $\cK$  be a
separable Hilbert space. Fix $\cB  =  \{ \varphi_n \}_{ n \in \N }\,$ an \onb  of $\cK$. 
From Eq. (\ref{frame}), one can deduce that
there exists a unique $T_{\F, \cB}  \in L(\cK , \H )$ such that
$T_{\F, \cB}  ( \varphi_n ) = f_n  \,$ for every  $ n \in \N$. 
We shall say that  $T_{\F, \cB}$ is a \textsl{preframe operator} for $\F$. 
Another consequence of Eq. (\ref{frame})
is that  $T_{\F, \cB} $ is surjective.
If one takes  the cannonical basis $\cE$ of $ \ell^2 (\NN )$, then 
$T_\F = T_{\F, \cE}$ is called the \textsl{synthesis operator} for $\F$.
\EOE
\end{fed}
\begin{rem} \label{cosas} \rm
Let $\F= \fram$ be a frame in $\hil$ and $T_{\F, \cB}\in \lkh$ a preframe operator
for $\F$, with $\cB =  \{ \varphi_n \}_{ n \in \N }\,$. Then  $T_{\F, \cB} ^*\in \lhk$  is   given by 
$\displaystyle T_{\F, \cB} ^*(x  ) = \sum_{n\in\N} \api x  , f_n \cpi \varphi_n\,$, for 
$x  \in \H$. It is 
an \textsl{analysis operator} for $\F$.
The operator $S_{\F}  = T_{\F, \cB} T_{\F, \cB} ^* \in \posop$,
called the \textsl{frame operator} of $\F$,  satisfies
$ 
S_\F  f=\sum_{n\in \N} \pint{f }{ f_n}f_n  \,$, for $f\in \H$. 
It follows from (\ref{frame}) that  $A_\F \, I \leq S_\F  \leq B_\F\,  I\, $.
So that $S_\F  \in \glh^+$.
Note that the frame operator $S_\cF$ does not depend on the preframe operator  chosen.
\EOE
\end{rem}

\begin{pro}\label{cotas optimas} \rm
Let $\F= \{f_j\}_{j\in J}$ be a frame sequence in $\hil$. 
Then the optimal frame constants for $\F$  are $A_\F =\gamma(T_\F)^2$ and $B_\F =\|T_\F\|^2$. \QED
\end{pro}

\begin{fed}\rm\label{intrinsecal}
Let $\F= \fram$ be a frame in $\hil$. The cardinal number
$\exe {\F } = \dim \ker T_\F $ 
is called the \textsl{excess} of the frame.
Holub \cite{[Ho]} and Balan, Casazza, Heil and Landau \cite{[BCHL]} proved that
\beq\label{balan}
\exe{\F }  = \sup \ \big\{ \ |I| :  I \inc \N \ \hbox { and } \
\{ f_n\}_{n\notin I} \ \hbox{ is still a frame for } \H \big\} \ .
\end{equation}
This characterization justifies the name ``excess of $\F$".
For every preframe operator $T_{\F, \cB} \in \lkh$ of $\F$, it holds that  
$\exe{\F } = \dim \ker T_{\F, \cB}\,$. 
The frame $\F$ is called a \textsl{Riesz basis} (or exact) 
if $ \exe{\F } = 0$, i.e., if
the preframe operators of $\F$ are invertible.
\EOE
\end{fed}

\section{Frames of subspaces, or fusion frames}
Throughout this section, $\hil$ shall be a fixed separable Hilbert space, 
and $I\inc \NN$ a fixed index set ($I = \NN$ or $I=\IN{n} := \{1, \dots , n \}$ for $n \in \NN$). 
Recall that $\linm$ denotes the space of bounded sequences of (strictly) positive numbers,  
which will be considered as weights in the sequel. The element $e \in \linm$ is the 
sequence with all its entries equal to $1$.

\subsection*{Prelimiaries}
Following Casazza and Kutyniok  \cite{[CasKu]}, we define:

\begin{fed} \rm
Let  $\ssec$ be a sequence    of closed subspaces of $\hil$, and let 
$w=  \{w_i\}_{i\in I} \in \linm$.  
\ben
\item
We say that $\sfram$ is a {\sl Bessel sequence}  of subspaces  (BSS) if there exists $B>0 $ such that  
\begin{equation}\label{BS de sub}
\sum_{i\in I}w_i^2 \, \|P_{W_i}f\|^2 \leq B\|f\|^2 \peso {for every $f\in \hil$ .} 
\end{equation}
where each $P_{W_i}\in \op$ is the orthogonal projection onto  $W_i\,$.
\item 
We say that 
$\cW_w$ is a {\sl frame of subspaces (or a fusion frame) for $\hil$}, 
and write that  $\cW_w$ is a FS (resp. 
FS for $\cS \subcer \hil$)  if there exist $A,B>0$  such that   
\begin{equation}\label{frame de sub}
A\|f\|^2\leq \sum_{i\in I}w_i^2 \, \|P_{W_i}f\|^2 \leq B\|f\|^2 \peso {for every $f\in \hil$ 
(resp. $f \in \cS$)\ ,} 
\end{equation}
The sharp constants for \eqref{frame de sub} are denoted by $A_{\cW_w} $ and $B_{\cW_w} \,$. 
\item  $\cW$  is a $minimal$ sequence if 
\beq\label{minimal}
 W_i \cap \gen{W_j : j \neq i} = \trivial  \peso{for every} i\in I \ .
 \end{equation} 
   \een
Suppose that $\cW_w$ is a fusion frame for $\hil$. Then 
\ben
\item [4.]
 $\cW_w$ is  a {\sl tight} frame if $A_{\cW_w} = B_{\cW_w} \,$, 
and {\sl Parseval} frame if  $A_{\cW_w} = B_{\cW_w} =1$.
\item [5.]
 $\cW_w$ is an {\sl orthonormal basis of subspaces} (shortly OBS) if 
 $ w=e$ and $W_i  \perp W_j$ for $i\neq j$.

\item [6.]  $\cW_w$ is {\sl Riesz basis of subspaces}  (shortly RBS) if 
$\cW $ is a minimal sequence. \EOE 
\een
\end{fed}

\medskip
\noi
The notions of synthesis, analysis and frame operators can be defined for BSS. 
But with a different structure of the Hilbert space of frame sequences, which now relies 
strongly in the sequence of subspaces $\ssec\,$.

\begin{fed}\label{defis} \rm
Let $\sfram$ be a BSS for $\hil$. 
Define the Hilbert space 
$$
\displaystyle \cK_\cW =\bigoplus_{i\in I}W_i \peso{with the $\ell^2$ norm}
\|g\|^2=\sum_{i\in I} \|g_i\|^2 
\ \ , \peso{for $g= (g_i)_{i\in I}\in \cK_\cW$ .}
$$ 
The {\bf Synthesis operator}: $\tfs \in L( \cK_\cW\, ,  \hil)$ is  defined by 
\[\tfs (g 
)=\sum_{i\in I} w_i \, g_i
\ \ , \peso{for $g = (g_i)_{i\in I}\in \cK_\cW$ .}\]
Its adjoint $\tfs^*\in L(  \hil, \cK_\cW )$ is called the {\bf Analysis operator} of $\FS\,$. It is easy to see that
$\tfs^*(f)=\{ w_i \, P_{W_i}f\}_{i\in I}\,$, for $f \in \hil$.
The {\bf Frame operator}: $\frs = \tfs \, \tfs^* \in \op^+$ satisfies the formula 
$ \frs f=\sum_{i\in I} w_i^2 \, P_{W_i}f$, for $f \in \hil$.
\EOE
\end{fed}

\begin{rem}\label{vectores} \rm
Let   $\ssec$  be  a sequence of closed subspaces of $\hil$, and let $w
\in \linm$. In  \cite{[CasKu]} the following  results were proved:
\ben
\item $\sfram $ is a BSS \iiff
the synthesis operator $\tfs $ is well defined and bounded. In this case, 
$$
 \peso{$\cW_w$ is a FS for $\hil$ (resp. for $\cS \subcer \hil$)} 
 \iff \peso{$\tfs$ is onto  (resp $R(\tfs ) = \cS$) \ .} 
$$
This is also equivalent to the fact that $\tfs^* $  is bounded from below. 
\een
If $\cW_w$ is a FS for $\hil$, then 
\ben
\item [2.]

$A_{\cW_w} = \gamma(\tfs )^2$ and $ B_{\cW_w} = \|\tfs\|^2$. 
So that $A_{\cW_w}  \cdot I\leq \frs \leq B_{\cW_w}  \cdot I$. 
\item  [3.]
$\cW_w$ is a RBS \iiff $\tfs $ is invertible (i.e. injective)
and $\cW_w$ is an OBS \iiff $w= e$ and $\tfs^* \tfs = I_{\cK_\cW} \, $. 
\item [4.]
$\cW_w$ is tight \iiff $\tfs \tfs^* = A_{\cW_w} \cdot I_\hil \, $, 
and $\cW_w$ is Parseval \iiff $\tfs $ is an coisometry (i.e. 
$\tfs \tfs^* = I_\hil \, $). 
\EOE
\een
\end{rem}

\medskip
\noi
We state another useful  result proved in \cite{[CasKu]}, which 
determines a relationship between frames of subspaces and frames of vectors.

\begin{teo}\label{CK} \rm
Let $\ssec$ be a sequence  of closed subspaces of $\hil$ and let $w \in \linm$. 
For each $i \in I$, let $\cG_i = \{f_{ij}\}_{j\in J_i}$ be a frame for  $W_i \,$. 
Suppose that 
$$
0<A=\inf_{i \in  I}   A_ {\cG_i} \peso {and}  B=\sup_{i \in  I} B_ {\cG_i} <\infty  \ .
$$ 
Let  $\cE_i = \{e_{ik}\}_{k\in K_i}$ be and \onb  for each  $W_i\,$. 
Then the following conditions are equivalent.
\ben 
\item  $\F = \{w_if_{ij}\}_{i\in I, j\in J_i}= \{w_i \, \cG_i \}_{i\in I} \, $ is a frame for $\hil$.
\item  $\cE = \{w_i e_{ik}\}_{i\in I, k\in K_i}= \{w_i \, \cE_i \}_{i\in I} \, $ is a frame for $\hil$.
\item  $\sfram $ is a frame of subspaces  for $\hil$. 
\een 
In this case, the bounds of $\cW_w$ satisfy the inequalities 
\begin{equation}\label{cottas}
\frac{\ds A_\F }{\ds B} \ \le A_{\cW_w }  = A_\cE   \peso
{and}
B_\cE  = B_{\cW_w } \ \le \ \frac{\ds B_\F }{\ds A}  \ .
\end{equation}
Also   $T_\cE = \tfs \,$, using the \onb 
$\cB = \{e_{ik}\}_{i\in I, k\in K_i}$ of $\cK_\cW =\bigoplus_{i\in I}W_i \,$. 
\QED
\end{teo}

\subsection*{Operators and frames}
Our next purpose is to characterize frames of subspaces as images of OBS under 
an epimorphism with certain properties. 

\begin{fed} \rm
Let $\sfram$ be a BSS for $\hil$, with synthesis operator $\tfs \,$. The excess of $\cW_w$ is
defined as:
$\exe{\cW _w } = \dim N(  \tfs ) \ . $
\end{fed}

\begin{teo}\label{cotas}
Let $\{E_i\}_{i\in I}$ be an OBS of $\cK$ and let $T \in L( \cK ,  \hil)$ 
be surjective. Suppose that 
$0< \ds \inf_{i \in I}\frac{\gamma(TP_{E_i})}{\|TP_{E_i}\|}\ $. Let 
$0<A, B<\infty$ be  such that, 
\begin{equation}\label{cota gama}
\frac{A}{B}\leq \frac{\gamma(TP_{E_i})^2}{\|TP_{E_i}\|^2} 
\peso{i.e. ,}  \frac{\|TP_{E_i}\|^2}{B}\leq \frac{\gamma(TP_{E_i})^2}{A} \ , \quad 
\forall \ i\in I \ .
\end{equation}
Denote $W_i = T(E_i)\subcer \hil$, for $i\in I$.  Let $w = \{w_i\}_{i\in I} \in \linm$ such that 
\begin{equation}\label{cota pesos}
\frac{\|TP_{E_i}\|^2}{B}\leq w_i^2\leq \frac{\gamma(TP_{E_i})^2}{A}
\peso{for each $i\in I$ .}
\end{equation}
Then the following statements hold:
\ben
\item The sequence  $\sfram $ is a FS  for $\hil$. 
\item Moreover, $\cW_w $ has bounds 
\begin{equation}\label{cota tutti}
\frac{\ds \gamma (T)^2 }{\ds B} \ \le A_{\cW_w }   \peso
{and}
B_{\cW_w } \ \le \ \frac{\ds \|T\|^2 }{\ds A}  \ .
\end{equation}
\item If $\ker T \cap E_i = \trivial$ for every $i\in I$, then  $\exe{\cW _w } = \dim \ker T 
\,$. 
\een
\end{teo}

\proof
Suppose  that \eqref{cota gama} and \eqref{cota pesos} hold for every $i\in I$. 
\ben 
\item Since $\gamma(TP_{E_i})>0$, then $W_i=TE_i$ is closed for every $i\in I$. 
Let $\{b_{ij}\}_{j\in J_i}$ 
be an \onb for each $E_i\,$. By Proposition \ref{cotas optimas},  
Eq. \eqref{cota gama} and Eq. \eqref{cota pesos}, every sequence 
 $\cG_i = \{w_i\inv \, T\, b_{ij}\}_{j\in J_i}$ is a frame for $W_i\,$ with 
 $$
 A_{\cG_i} = w_i ^{-2}  \gamma(TE_i)^2 \ge A 
 \peso {and}
 B_{\cG_i} = w_i ^{-2}  \|TE_i \|^2 \le B \ . 
 $$
 On the other hand, since 
$\{b_{ij}\}_{i\in I, j\in J_i}$ is a \onb for $\cK$, and $T$ an epimorphism, the sequence 
 $\F = \{Tb_{ij}\}_{i\in I, j\in J_i}$ is a frame for $\hil$.
Finally, since $\F=\{w_i(w_i^{-1}Tb_{ij})\}_{i\in I, j\in J_i} = \{w_i \, \cG_i \}_{i\in I} \, $, 
Theorem \ref{CK} implies that $\cW_w $ is a FS for $\hil$. 
\item 
Eq. \eqref{cota tutti} follows from Eq. 
\eqref{cottas} and the fact that 
$A_\F = \gamma(T)^2 $ and $B_\F = \|T\|^2$. 
\item Suppose that $\ker T \cap E_i = \trivial$ for every $i\in I$. 
Then $\ker TP_{E_i} = E_i \orto$ and $\gamma (TP_{E_i} ) \ \|z\| \le \|TP_{E_i} z \|$ 
for every $z \in E_i\, $. By Eq. \eqref{cota pesos},  for every $x \in \cK$ and   $i\in I$, 
$$
A\rai w_i \ \|P_{E_i} \ x \| \le \gamma (TP_{E_i} ) \ \|P_{E_i} \ x  \| \le \|T P_{E_i} \ x  \| 
\le B\rai w_i \|P_{E_i} \ x \|\ ,
$$
and $\|x\|^2 =  \sum_{i\in I} \|P_{E_i} x\|^2$. 
Let $\cK_\cW = \bigoplus_{i \in I } W_i \,$ (the domain of $\tfs \,$).
Observe that  $T(E_i ) = W_i $ for every $i \in I$.  
Therefore the map 
$$
V : \cK \to \cK_{\cW}  \peso{given by} V x  
= \big(\ w_i \inv \ T ( P_{E_i} \ x ) \,\big)_{i \in I}\ \ , \peso{for} x \in \cK \ ,
$$
is well defined, bounded and invertible. By the definition of 
the synthesis operator $\tfs\,$, and the fact that $x = \sum _{i \in I} P_{E_i} x $, 
for every $x \in \cK$, we can deduce that $\tfs \circ V  = T$. 
Therefore $\dim \ker T = \dim V\inv (\ker \tfs )= \dim \ker \tfs = \exe{\FS}$. 
\QED
\een

\medskip
\noi
Example  \ref{gama no frame} shows  a surjective operator $T$ and an 
OBS $\cE = \{E_i \}_{i\in I}$ such that  $\gamma(TP_{E_i})>0$ for every $i\in I$,  
but the sequence $\cW_w = (w, \cW) $ fails to be a FS for every $w\in \linm \,$. Hence 
 $T$ and $\cE$ do not satisfy Eq.  \eqref{cota gama}.

However, Eq.  \eqref{cota gama}  is not a necessary  condition in order 
to assure  that $\pesos \neq \vacio $ (see Definition \ref{pesos}), if $\cW = T\cE$. 
In Example \ref{no dos} we show a FS wich is the image of an 
OBS under an epimorphism which doesn't satisfy Eq. \eqref{cota gama}.

\medskip

\begin{rem}\label{wi}
If $\sfram$ is a FS for $\hil$, then its synthesis 
operator $\tfs \,$, defined as in Definition \ref{defis} clearly satisfies  
Eq. \eqref{cota gama}. Moreover, it holds that 
$$
\tfs g = w_i g  \peso{for every $g \in E_i \,$, the copy of $W_i$ in $\cK_\cW$\ .}
$$ 
Hence  $ \gamma (\tfs P_{E_i} ) = \|\tfs P_{E_i} \|=  w_i $ for every $i\in I$.  
\EOE
\end{rem}

\medskip

\begin{rem} \rm Let $\sfram$ be a FS for $\hil$, and let $G \in \glh$. 
In \cite{[CasKu]}, \cite[Thm 2.11]{[CasKuLi]} and \cite [Thm 2.4]{[Gav]} it is proved 
that  $G\cW_w = (w_i \, ,  GW_i ) _{i \in I}$ must be also
 a FS for $\hil$. 
We give a short proof of this fact, including extra information about the bounds and the 
excess of $G\cW_w \,$, in order to illustrate the techniques given by Theorem \ref{cotas}. \EOE
\end{rem}
\begin{cor}\label{conG}
Let $\sfram$ be a FS for $\hil$, and let $G\in L(\hil , \hil_1 ) $ be {\bf invertible}. Then 
$G\cW_w = (w_i \, ,  GW_i ) _{i \in I}$ is a FS for $\hil_1\,$, 
which satisfies that   $\exe{\cW _w } = \exe{G \, \FS}$,  
\[
\left( \|G\| \, \|G\inv\|\right)^{-2} \ A_{\cW_w } \ \le A_{G\cW_w }  \peso{and} 
B_{G\cW_w } \ \le \ \left( \|G\| \, \|G\inv\|\right) ^2 \ B_{\cW_w }\ .
\] 
\end{cor}
\proof 
Denote by $E_i$ 
the copy of each $W_i$ in $\cK_\cW= \bigoplus_{i \in I } W_i \,$.  
Define $ T = G \, T_{\cW_w}  \in L(\cK_\cW \, , \hil_1)$, which is clearly 
surjective (since $T_{\cW_w} $ is). 
By  Eq. \eqref{desi1} and Remark \ref{wi}, 
$$
\gamma (T  P_{E_i} ) \ge \gamma(G)  \cdot \gamma (T_{\cW_w}  P_{E_i} )  = \gamma(G)  \ w_i  \peso{and}
\|T  P_{E_i} \| \le \|G\| \, \|T_{\cW_w}  P_{E_i} \| = \|G\| \ w_i  \ ,
$$ 
for every $i \in I$. In particular,  $T (E_i ) \subcer \hil_1\,$. Then, we can apply Theorem \ref{cotas}
for $T$ with constants $A = \gamma(G)^2$ and $B = \|G\|^2$. Indeed, for every $i \in I$, we have seen that
$$
\frac{\ds \gamma(G)^2 }{\ds \|G\|^2} \le \frac {\ds \gamma (T P_{E_i} )^2}{\ds \|T  P_{E_i} \|^2}
\peso{and} 
 \frac{\ds \|T  P_{E_i} \|^2}{\ds \|G\|^2} \ \le \  w_i ^2 \ \le 
 \ \frac{\ds \gamma (T P_{E_i} )^2}{\ds \gamma(G)^2 } \ .
$$
Therefore, $G\cW_w = (w_i \, ,  GW_i ) _{i \in I}$ is a FS for $\hil_1$ by Theorem \ref{cotas}.
In order to prove the bound inequalities, by Eq. \eqref{desi1} and item 2 of Remark \ref{vectores}  we have that 
$$
\gamma (G\tfs ) \ge \gamma(G) \, \gamma(\tfs ) = \|G\inv\|\inv \, A_{\cW_w}\rai   \peso{and} 
\|G\tfs \| \le \|G\| \, \|\tfs\| = \|G\| \, B_{\cW_w}\rai  \ .
$$
Now apply Eq. \eqref{cota tutti}   of Theorem \ref{cotas}
with our constants $A = \|G\inv\|^{-2} $ and $B = \|G\|^2$.
It is easy to see that $\ker T = \ker \tfs\,$. Then $\ker T \cap E_i = \trivial$ ($i \in I$). 
By Theorem \ref{cotas}, we deduce that $\exe{\FS}= 
\dim \ker \tfs = \dim \ker T = \exe{G \, \FS}$. 
\QED

\section{Admissible weights}

\begin{fed}\label{pesos}   \rm
We say that  $\ssec$ is a {\sl generating} sequence  of $\hil$, if 
$W_i \subcer \hil$ for every $i \in I$, and  $\gen {W_i  : i\in I} = \hil$.
In this case, we define
\[
\pesos= \big\{ \ w \in \linm : \,  \sfram \text { is a FS for } \hil \big\} 
\inc \linm  \ ,
\]
the set of $admisible$ sequences of weights for $\cW$. \EOE
\end{fed}

\medskip
\noi
It is apparent that, if $\sfram$ is  a FS for $\H$, then 
$\ssec$ is a generating sequence. Nevertheless, 
in Examples \ref{gama no frame}  and \ref{sin int} we shall see that there exist generating 
sequences $\cW= \{W_i\}_{i\in I}$ for $\hil$ such that $\pesos = \vacio$. Recall that we denote by 
\beq\label{bsi2}
\bsi = \{ \, \{w_i\}_{i\in I} \in \linm : \inf _{i \in I} w_i > 0 \} = \linm \cap Gl (\ell^\infty(I)\,) \ .
\end{equation}

\begin {pro} \label{los pesos} 
Let $\cW = \{W_i\}_{i\in I}$ be a generating sequence of $\hil$. 
\ben
\item If $w \in \pesos $,  then \ $a  \, w \in \pesos $ and 
$ \exe{\cW _w } = \exe{\cW _{a\, w} }\,$, for every $a \in \bsi $.
\item If $\cW_w = (w , \cW)$ is a RBS, for some $w \in \linm$, then
$ \pesos =  \bsi$, and $(a, \cW)$ is still a   RBS for every $a \in  \bsi$. 
In particular, $(e, \cW )$ is a RBS. 
\item Let $G \in \glh$. Then  $\pesos=\cP\left(\{GW_i\}_{i\in I}\right)$. In other words, a sequence 
$w \in \linm $ is admisible for $\cW$ \iiff it is admisible for $G\cW$. 
\een 
\end{pro}

\proof Let $\cK_\cW = \bigoplus_{i \in I } W_i \,$, and denote by $E_i\subcer \cK_\cW\, $ 
the copy of each $W_i$ in $\cK$. 
\ben
\item 
For every  $a \in  \bsi$, consider the SOT limit $D_a = \ds \sum_{i\in I} a_i P_ {E_i} \, $. 
Then $D_a \in \glkw ^+$. Therefore, if $\tfs \in L(\cK_\cW \, , \hil)$ is the synthesis 
operator of $\FS$, then $\tfs \circ D_a $ is, by definition, 
 the synthesis operator of $(a \, w, \cW)$. Since $\tfs \, D_a $ is bounded and  surjective, then  
 $(a \, w, \cW)$ is also a FS. Note 
 that $N(  T_{\cW _{a\, w} } )= N(  \tfs \, D_a ) =D_a\inv (N(  \tfs )\, )  $. 

\item If $\FS$ is a RBS for $\H$, then $\tfs$ is invertible. Since 
$\tfs x = w_i x$ for $x \in E_i\,$, then  $w_i \ge \gamma (\tfs ) = A_{\FS}\rai $
for every $i \in I$. 
This implies that $w \in  \bsi$. 
Observe that $w  \bsi=  \bsi$ (because $w\inv \in  \bsi$). Then $ \bsi \inc \pesos$ by item (1). 
But, for every  $a \in \pesos$, we have that $\cW_a$ is a RBS, because $\cW$ is still minimal. 
Then $s\in \bsi$. 
\item  Apply Corollary \ref{conG} for $G$ and $G\inv$. \QED
\een

\begin{fed}\rm 
Let $\cW = \{W_i\}_{i\in I}$ be a generating sequence of $\hil$. Given $v, w \in \pesos$, 
we say that $v$ and $w$ are $equivalent$ if there exists $a \in \bsi $ such that $v = a \cdot w$. \EOE
\end{fed}
\begin{rems} \label{muchos pesos} 
Let $\cW = \{W_i\}_{i\in I}$ be a generating sequence of $\hil$. 
\ben
\item 
By Proposition \ref{los pesos}, if  $w\in \pesos$, then its whole equivalence class $w \cdot \bsi \inc \pesos$. 
\item On the other hand, in Example \ref{frame no gama} below we shall see that there exist 
generating sequences $\cW $ of $\hil$  with infinite not equivalent sequences 
$w \in \pesos$. 
%
\item If  $\cW_w$ is a RBS for $\hil$, then by Proposition \ref{los pesos}  
all the admissible sequences for $\cW$ are equivalent to $w$, since  
$ \pesos =  \bsi\,$. 
Since $\cW_v = (v , \cW)$ is a RBS for $\hil$ for every $v \in  \bsi$,   
 from now on we  will not mention the weights. We just say that the sequence of subspaces
$\cW$ is a Riesz basis of subspaces.  
\item By definition, if $\cW$ is a RBS, then it is a minimal sequence. 
 Nevertheless, in Example
 \ref{sin int}, we shall see that there exist minimal sequences which are generating for $\hil$, 
 but with $\pesos = \vacio$. 
\EOE
\een
\end{rems}

\begin{pro} \label{onRi} Let $\cE = \{E_i \}_{i\in I} $ be a OBS for  $\hil$.
Let $G \in L( \hil ,  \hil_1 )$ be an invertible operator. Then the sequence $\cW = (GE_i) _{i \in  I}$
is a RBS for $\hil_1\,$.
\end{pro}
\proof It is a consequence of Corollary \ref{conG}. \QED

\begin{rem}\label{no ti} 
It is well known (and easy to verify) that for a frame $\cF=\{f_i\}_{i\in I}$ in $\hil$, the 
sequence  $\{S_{\cF}^{-1/2}f_i\}_{i\in I}$ is a Parseval frame. 
Nevertheless, if $\sfram \,$ is a FS, then $\frs\mrai \cW_w\,$ 
may be not a Parseval FS (see Example \ref{frame no gama} below), 
neither allowing to change the sequence of weights. 
Even worse, there  exist frames of subspaces $\FS = (w , \cW)$ for $\hil$ such that  
the sequence $(v, G \, \cW)\,$ fails to be a Parseval FS for $\hil$ for every $G \in \glh$ 
and $v \in \linm$  (see Example  \ref{otros2}). In the next Proposition 
we shall see that the situation is different for a RBS of $\H$: 
\EOE
\end{rem}

\begin{pro}\label{en RBS funca}
Let $\ssec$ be a RBS for $\hil$. Then, for every $w \in \bsi$, 
the sequence  $\{ S_{\FS}^{-1/2}W_i\}_{i\in I}\,$ is an \onb  of subspaces.
\end{pro}
\proof
Let $\{e_{ik}\}_{k\in K_i}\,$ be an \onb of each $W_i\,$. 
According Theorem \ref{CK}, the sequence  $\cE = \{w_i e_{ik}\}_{i\in I, k\in K_i}$ 
is a Riesz basis  of $\hil$ and $T_{\cE}=T_{\FS}\,$.  
Hence  the sequence  $\{w_iS_{\cE}^{-1/2} e_{ik}\}_{i\in I, \, k\in J_i}\, $ is an 
\onb  for $\H$. Since $S_{\FS}=S_{\cE}\, $ and $\{w_i S_{\FS}^{-1/2} e_{ik}\}_{k\in K_i}\,$
is a \onb of each subspace $S_{\FS}^{-1/2}W_i\,$,  then 
$\{S_{\FS}^{-1/2}W_i\}_{i\in I}\,$ is an OBS for $\H$.
\QED

\section{Projections and frames}

In this section we obtain a generalization of  two results of \cite{[ACRS2]}, which relates FS 
(including the computation of their weights) and oblique projections
(see also \cite{[AsKho]}). Unlike for vector frames, all the results are in ``one direction''. 
The converses fail in general 
(see Example \ref{proyeccion no frame} and Remarks \ref{converse1} and \ref{converse2}).

\begin{teo}\label{con R}
Let $\sfram $ be a  FS for $\hil$. Then there exists a Hilbert space 
$\cV \supseteq \hil$ and a {\bf Riesz  basis} of subspaces $ \{B_i\}_{i\in I}$ for  $\cV$ such that
\[ 
P_{\hil}  (B_i) \ = \ W_i  
\peso{and} 
A_{\cW_w}\rai \  \|P_{\hil}P_{B_i}\|\leq w_i\leq B_{\cW_w}\rai \ \|P_{\hil}P_{B_i}\| 
\peso{\rm for every $i\in I$ .}
\] 
This means that the new sequence of weights $v_i = \|P_{\hil}P_{B_i}\|$, $i \in I$, is equivalent to $w$.
Also, we can compute $\exe{\FS } = \dim \cV\ominus \H$.
\end{teo}
\proof
Denote by $E_i$ 
the copy of each $W_i$ in $\cK_\cW\,  = \bigoplus_{i \in I } W_i \,$. 
Let $\tfs \in L(\cK_\cW\, ,  \hil)$ be the synthesis operator for $\cW_w\,$. 
Denote by  $ \cN = N(\tfs)$ and 
$\cV=\hil\oplus \cN $. We can identify $\hil$ with $\hil\oplus \{0\} \subcer \cV$. 
Let |
\beq \label{el U}
U\ :\cK_\cW\,  \ra \cV \peso{given by} 
U(x)=\tfs x\oplus  \gamma (\tfs )  \ P_{\cN}\ x   \quad , \quad x \in \cK_\cW \ .
\end{equation}
Since $\cK_\cW\, =\cN \orto \perp  \cN$ and $\tfs \big|_{\cN\orto} : 
 \cN\orto \to \hil$ is invertible, we can  deduce that  $U$ is bounded and invertible. 
Moreover, it is easy to see that 
\beq\label{el U 2}
\|U^{-1}\|^{-1}=\gamma (U)=\gamma(\tfs ) = A_{\cW_w}\rai \peso{and} \|U\|=\|\tfs\| = B_{\cW_w}\rai\ .
\end{equation}
By Proposition \ref{onRi}, the sequence $\{B_i\}_{i\in I} = 
\{U(E_i)\}_{i\in I}$ is a  RBS for $\cV$.
Observe that 
\[
P_\H (B_i ) = P_\H 
U(E_i) =\tfs (E_i)\oplus\{0\}=W_i\oplus \{0\} \sim W_i \ \ , 
\peso{\rm  for every $i\in I$  .}
\]
Let $y$ be an unit vector of $B_i=U (E_i)$. Then $y = Ux$ with 
$x\in E_i\,$. We have that
\[\gamma(U)\|x\|\leq \|Ux \| = \|y\| = 1 \leq \|U\| \, \|x\|
\ .
\]
Recall that $E_i$ 
is the copy of $W_i$ in $\cK$.  If $x \in E_i\,$,  we denote by $x_i$ its component in $W_i\,$
(the others are zero). 
Using that $\|P_{\hil} \, y\|= \|\tfs \, x \| = w_i\|x_i\|=w_i\|x\|$ and Eq. \eqref{el U 2}, 
we can conclude that 
for every such $y$ (i.e. any unit vector of $B_i$), 
$$
 A_{\cW_w}\rai\  \|P_{\hil} \ y \| \ = \ \gamma(\tfs )  \,  \|P_{\hil} \ y \| =   
 w_i \ \gamma (U) \,  \|  x \| \le w_i 
 \ \implies \  A_{\cW_w}\rai\ \| P_{\hil} P_{B_i} \| \le w_i \ .
 $$
 Similarly, 
$
w_i \le w_i \ \|U\|\,  \| x \| = B_{\cW_w}\rai\  \|P_{\hil} \ y \| \le B_{\cW_w}\rai\ 
\| P_{\hil} P_{B_i} \|  
$. 
\QED

\bigskip
\noi
As a particular case of Theorem \ref{con R}, we get a result proved by Asgari and Khosravi 
\cite{[AsKho]} (see also \cite{[CasKuLi]}), with some information extra:

\begin{cor}\label{con ON}
Let $\sfram $ be a Parseval FS for $\hil$. Then there exists a Hilbert space 
$\cV \supseteq \hil$ and an {\bf orthonormal basis} of subspaces 
$\{F_i\}_{i\in I}$ for  $\cV$ such that
\[ 
P_{\hil} ( F_i) \ = \ W_i \peso{and} w_i = \angf{\hil}{F_i} = \|P_\hil \, P_{F_i} \| 
\peso{\rm for every $i\in I$ .}
\] 
\end{cor}
\proof
We use  the notations of the proof of Theorem  \ref{con R}. If $\cW_w$ is 
Parseval, then 
$A_{\cW_w} = B_{\cW_w} = 1$. By Eq. \eqref{el U 2}, this implies that the operator 
$U\in L(\cK  ,  \cV)  $ defined in Eq. \eqref{el U} becomes unitary (it is an invertible isometry). 
Hence, in this case,  the sequence $\{F_i\}_{i\in I} = \{U(E_i)\}_{i\in I}$ is a  \onb of 
subspaces  for $\cV$. Also, by Theorem  \ref{con R}, we have that 
$w_i  = \|P_\hil \, P_{F_i} \| $ for every $i \in I$. 
It is easy to see that  $F_i \cap (\hil \oplus \trivial ) \neq \trivial $ implies that 
$w_i = 1$ and $F_i \inc (\hil \oplus \trivial )$   (because $U$ is unitary). 
Then, we can deduce that    $\|P_\hil \, P_{F_i} \|
= \angf{\hil}{F_i} $ for every $i \in I$. \QED

\begin{rem} \label {converse1} \rm 
Although the converse of Corollary \ref{con ON}   fails in general, it 
holds with some special assumptions, based on Theorem \ref{cotas}: 
If $\cE =  \{E_i\}_{i\in I}$ is a 
OBS   for  $\cV\supseteq \hil $ such that
$\ds  0 < \inf_{i \in I}  \frac {\gamma (P_{\hil}P_{E_i})}{\|P_{\hil}P_{E_i}\|} \, , 
$ then $\pesos \neq \vacio $, where $W_i = P_{\hil}  (E_i)$, $i \in I$. 
Moreover, as  in Theorem \ref{cotas}, it can be found a concrete  $w \in \pesos$. 
Nevertheless, we can not assure that $\cW_w$ 
is a Parseval FS.  
\EOE
\end{rem}

\medskip
\noi
The following theorem  is closely related with a result proved by Casazza, Kutyniok and Li in 
\cite[Thm. 3.1]{[CasKuLi]}. 

\begin{teo}\label{T:poyeccionoblicua a frames}
Let $\sfram$ be a FS for $\hil$ such that $1\leq A_{\cW_w}\,$. 
Denote by $\cV=\hil \oplus \cK_\cW\, $.
Then there exist an
{\bf oblique} projection $Q \in L(\cV ) $ with $R(Q) = \H \oplus \{0\}$
and an {\bf orthonormal system}  of subspaces $ \{B_i\}_{i \in I}$ in $\cV$,  such that
$$
W_i \oplus 0 =Q  (B_i) \peso{and}
w_i = \|Q \, P _{B_i} \| = \gamma (Q \, P _{B_i} ) \peso{for every \ $i \in I$ .}
$$
Moreover, if $\exe {\cW_w} =   \infty$,  then  the sequence $\{B_i\}_{i \in I}$
can be supposed to be an {\bf orthonormal basis} of subspaces of $\cV$. 
\end{teo}

\begin{proof}
Write $\tfs = T$. 
By hypothesis, $TT^* = \frs \geq A_{\cW_w} I \ge I $. Denote by 
$$
X=(TT^*-I)^{1/2} \in \cam \ .
$$
Consider the (right) polar decomposition $T=|T^*|V\,$, where $V \in L( \cK_\cW \, , \hil) $ is a 
 partial isometry with initial space $N(T)^\bot$ and final space $\hil$, so that 
$V V^*=I_{\hil}\,$. 
Consider the ``ampliation" $\tilde T \in L(\cK_\cW \, , \cV) $ given  by
$\tilde T x=T x\oplus 0$.  Then 
$\tilde  T \tilde T^* = \bm {cc} TT^* & 0 \\ 0 & 0 \em \barr{l} \hil \\ \cK_\cW
\earr\in L(\cV )$. 
Define 
\[Q=\bm {cc} I_\hil & XV  \\ 0 & 0 \em \barr{l} \hil \\ \cK_\cW
\earr \in L(\cV ) \ . 
\]
Then it is clear that $Q$ is an oblique  projection with $R(Q) = \hil \oplus 0$. Moreover,
$$
QQ^*= \bm {cc} I_\hil + XX^* & 0 \\ 0 & 0 \em  =
\tilde  T \tilde T^*  \implies |Q^*|=|\tilde T ^*| \ .
$$
 Define $U\in L( \cK_\cW \, , \cV ) $  by
\beq\label{laU}
Ux=V P_{N(T)^\bot}x\oplus P_{N(T )}x \ \  , \peso{for $x\in \cK_\cW $  .}
\end{equation}
Then $U$ is an {\bf isometry}, because the initial space of $V$ is $N(T)^\bot$. 
Note that also $\tilde T =|\tilde T^*|U$. The partial
isometry of the right polar decomposition of $Q$ extends to an
unitary operator $W$ on $\cV$, because
$\dim N(Q) = \dim R(Q)^\bot $. Moreover, $Q=|Q^*|W$. Then
\[\tilde T =|\tilde T ^*|U=|Q^*|U=Q\ W^*U .\]
Therefore, if we consider the OBS $\{E_i\}_{i\in I}$ of $\cK_\cW\, $,

$$ W_i =T (E_i ) \sim T (E_i) \oplus 0  =
\tilde T (E_i)  =Q W^*U (E_i) =Q(B_i) \ , \quad i \in I \ , 
$$ 
where $\{B_i\}_{i \in  I} = \{W^*UE_i \}_{i \in  I}\, $, which is clearly an orthonormal 
system in $\cV$. 
If $y \in B_i$  is an unit vector, 
then $y = W^*U x$ for $x \in E_i\,$ with $\|x \| = 1$, and 
$$
w_i = \|Tx \| = \| QW^*U x \| = \|Q y \| \implies w_i = \|Q \, P _{B_i} \| = \gamma (Q \, P _{B_i} )  \ .
$$
Suppose now that $\dim N( T) = \infty$. Then 
the isometry $U$ defined in equation (\ref{laU})
can be changed to an unitary operator from $\cK_\cW$ onto $\cV$, still
satisfying that $\tilde T =|\tilde T^*|U$.
Indeed, take
$$ U'x= V P_{N(T )^\bot }x\oplus Y \ P_{N(T)} x \ , \peso{for} x \in \hil \ ,  
$$ 
where $Y \in L(\cK_\cW )$ is a partial isometry with
initial space $N(T)$ and final space $\cK_\cW\,$.
It is easy to see that $U'$ is unitary. 
Then the sequence $\{ B'_i\}_{i \in I} = \{W^*U'E_i\}_{i \in I}\,$  turns to be an OBS for $\cV$.
\end{proof}

\begin{rem} \label{converse2} \rm As in Remark \ref{converse1}, 
it holds a kind of converse for Theorem  \ref{T:poyeccionoblicua a frames}, 
i.e., if $ \ds \inf _{i\in I} \frac{\gamma (Q \, P _{B_i} )}{ \|Q \, P _{B_i} \|} >0$, 
then $\cP( \{ Q(B_i)\}_{i\in I} )\neq \vacio$.   \EOE
\end{rem}

\section{Refinements of frames of subspaces}
In \cite{[CasKu]} it is shown by an example that a FS  with $\exe{\FS }>0$ can be exact, i.e. 
$(w_i , W_i)_{i\in J}$ is not a  FS, for every proper $J\subset I$.
This situation is possible because the excess of the frame can be contained properly in some $W_i\in \FS$, so if we ``erase" any of the subspaces if $\FS$, this new sequence  is not generating anymore.

Then, the notion of ``excess" is not the same as for vector frames, in the sense 
of Definition \ref{intrinsecal} and Eq. \eqref {balan}. 
In this section, we introduce the notion of refinements of subspace sequences, which shall work as 
the natural way to recover the connection between excess and erasures. The results 
of this section are closely related with those of \cite [Section 4]{[CasKu2]}.

\begin{fed} \rm 
Let $\ssec$ be a sequence of closed subspaces. 
\ben
\item A {\sl refinement} of $\cW$ is a sequence $\refi$ of closed subspaces such that 
\ben
\item $J \inc I$.
\item $\trivial \neq V_i\inc W_i$ for every $i\in J$.
\een
\een
In this case we use the following notations:
\ben
\item [2.] The excess of $\cW$ over $\cV$ is the cardinal number 
$$
\exe{\cW, \cV} = \sum_{i\in \ J }\dim ( W_i\ominus V_i ) + \sum_{i\notin \ J}\dim  W_i  \ .
$$
\item [3.] If $w \in \pesos $, we say that $\cV_w = (w_i , V_i)_{i\in \ J}$ is 
a {\sl FS refinement} (FSR) of $\FS$ if $\cV_w $ is a FS for $\hil$. \EOE
\een
\end{fed}

\begin{rem} \label{trans}
It is easy to see that, if  $\cV$ is a refinement of $\cW$ and 
$\cV'$ is a refinement of $\cV$, then 
$\cV'$ is a refinement of $\cW$ 
and $\exe{\cW, \cV'} = \exe{\cW, \cV} +\exe{\cV, \cV'} $. \EOE
\end{rem}

\medskip
\noi
The next result uses basic Fredholm theory. We refer to J. B. Conway book 
\cite[Ch. XI]{[CW]}.

\begin{lem}\label{proy}
Let $\sfram$ be a FS for $\hil$ and let $\refi$ be a refinement of $\cW$. 
We consider  $\cK_\cV = \oplus_{i \in J} \cV_i\, $ as a 
subspace of $\oplus_{i \in I} \cW_i = \cK_\cW \,$. Then 
\ben
\item 
$\exe{\cW, \cV} = \dim \cK_\cV \orto = \dim \ker P_{\cK_\cV} \,$.
\item $\cV_w = (w_i , V_i)_{i\in \ J}$ 
a FS refinement of $\FS$ \sii $\tfs P_{\cK_\cV}$ is surjective. 
\een
In this case, we have that 
\ben
\item [3.] $\exe{\cW, \cV} \le \exe{\cW_w }$. 
\item [4.] If $\exe{\cW, \cV }<\infty$, then $\exe{\cV_w } = \exe{\cW_w }  -\exe{\cW, \cV}  $. 
\een
\end{lem}
\proof 
For each $i \in I$, denote by $E_i $ (resp. $F_i$) the copy of $\cW_i $ (resp. $\cV_i \,$, 
or $F_i = \trivial$ if $i \notin J$) in $\cK_\cW\,$. Then $\cK_\cV \orto = 
\oplus_{i \in I} E_i \ominus F_i \, $, showing (1). 
Denote by $P = P_{\cK_\cV}\, $. 
By construction, $T_{\cV_w} = \tfs \big|_{\cK_\cV} = \tfs \big|_{R(P)} \in L(\cK_\cV\, , \H)$. 
Then $R(\tfs P) = R(T_{\cV_w}) =\hil$ \sii $\cV_w$ 
a FS refinement of $\FS$. In this case,  $\trivial = \ker P \tfs^* \,$. Since $R(\tfs^* ) = 
\ker \tfs \, \orto$, then 
$$
\ker \tfs \, \orto \cap \ker P = \trivial \implies 
\exe{\cW, \cV} = \dim \ker P  \le  \dim \ker \tfs =  \exe{\cW_w } \ .
$$
Observe that $\tfs $ is a semi-Fredholm operator, with 
$\ind{\tfs }= \dim \ker \tfs -0 = \exe{\cW_w }$. 
If $\exe{\cW, \cV } <\infty$, then $P$ is a Fredholm operator, with $\ind{P}=0$. Hence, we have that  
$  \exe{\cW_w } = \ind{\tfs } + \ind{P} = \ind{\tfs P} =   \dim \ker\tfs P $. 
Finally, since $T_{\cV_w} = \tfs \big|_{\cK_\cV}\,$, 
$$
\exe{\cV_w } = \dim \ker T_{\cV_w} = \dim \ker\tfs P - \dim \ker  P = 
\exe{\cW_w }  -\exe{\cW, \cV}   \ ,
$$
which completes the proof. 
\QED

\begin{lem}\label{exari}
Let $\sfram$ be a FS for $\hil$ with 
$\exe{\cW_w }>0$. Then there exists  a FS refinement $\cV_w = (w_i , V_i)_{i\in \ J}$ 
 of $\FS$ with $\exe{\cW, \cV} =1$. 
\end{lem}
\proof
For each $i \in I$, denote by $E_i $ the copy of $\cW_i $  in $\cK_\cW\,$.
Suppose that there is  no FS refinement $\cV_w$ 
 of $\FS$ with $\exe{\cW, \cV} =1$. Then, by Lemma \ref{proy}, 
 for every $i \in I$ and every 
 unit vector $e \in E_i \,$, it holds that $R(\tfs P_{\{e\}\orto } ) \neq \hil$.  
 By Proposition \ref{propiedades elementales de los angulos} and 
Eq. \eqref{desi1},  
  $$
  \angf{N(\tfs )}{\{e\}^\bot } = \angf{N(\tfs )\orto}{\sspan\{e\}} <1 \implies R(\tfs P_{\{e\}\orto })\subcer\H  \ .
  $$
Take $x_e \in R(\tfs P_{\{e\}\orto })\orto = \ker P_{\{e\}\orto} \tfs^* $ an unit vector. Then 
$0\neq \tfs^* x_e \in \sspan\{e\}$, i.e., $e \in R(\tfs ^*)$. 
This implies that $\cup _{i\in I} E_i \inc R(\tfs ^*)$ (which is closed), so that 
$\tfs^*$ is surjective and $\exe {\FS } = 0$. \QED

\begin{teo} \label{raros0} Let $\sfram$ be a FS for $\hil$. Then \rm
\beq\label{eccesos}
\exe{\FS } = \sup \Big\{ \exe{\cW, \cV} \ :  \   \cV_w = (w_i , V_i)_{i\in \ J}
 \peso{is a FS refinement of} \cW_w  \ \Big\} \ .
\end{equation}
In particular, if  $E(\FS )=\infty$, then, for every $n\in \N$, there exists a FS refinement 
$\cV_w = (w_i , V_i)_{i\in \ J}$ of $\cW_w$ such that $\exe{\cW, \cV} = n$. 
\end{teo}
\proof Denote by $\alpha$ the supremum of Eq. \eqref{eccesos}. 
Observe that item 3 of Lemma \ref{proy} says that $\alpha \le \exe{\FS }$. If $\exe{\FS }<\infty$, 
combining Remark \ref{trans}, Lemma \ref{exari} and item 4 of Lemma \ref{proy}, one obtains an inductive argument which 
shows that $\alpha \ge \exe{\FS }$.  If $\exe{\FS }= \infty$, a similar
inductive argument shows that, 
 for every $n\in \N$, there exists a FS refinement 
$\cV_w $ of $\cW_w$ such that $\exe{\cW, \cV} = n$.  \QED

\begin{cor} \label{refis1}
Let $\sfram$ be a FS of $\hil$ such that $E(\cW_w )<\infty$. Then 
\ben
\item The sequence $w\in \bsi$.
\item There exists a FS refinement $\cV_w = (w_i , V_i)_{i\in \ J}$  of $\cW_w$ such that:
\begin{enumerate}
\item $\cV$ is a RBS for $\hil$.
\item $\exe{\cW , \cV} =E(\cW_w)$.
\end{enumerate}
\een
\end{cor}
\proof
By Theorem \ref{raros0}, there exists a FS refinement $\cV_w  = (w_i , V_i)_{i\in \ J}$  
of $\cW_w$ such that $\exe{\cW , \cV} =E(\cW_w)$. By item 4 of Lemma \ref{proy}, $E(\cV_w) =0$. 
This means that $\cV_w$ is a RBS for $\H$. Then, 
by Proposition \ref{los pesos},
the sequence $\{w_i\}_{i \in J } \in \linmJ^*$. Since $\exe{\cW , \cV} < \infty$, 
then   $I\setminus J$ is finite, and we get that also $w \in \bsi$. 
\QED

\begin{cor} \label{raros}
Let $\sfram$ be a FS for $\hil$ such that  $\exe {\cW_w} <\infty$. Then 
$$
\pesos =  \bsi  \peso{and} \exe {\cW_v } = \exe {\cW_w } 
\peso{for every other} v \in \pesos \ .
$$ 
\end{cor}
\proof 
By Corollary \ref{refis1}, we know that $w \in \bsi$. By Proposition \ref{los pesos}, we deduce that 
$\bsi \inc \pesos$. 
Let $\cV_w  = (w_i , V_i)_{i\in \ J}$ be a FS refinement of $\FS$ which is a RBS for $\H$, 
provided by Corollary \ref{refis1}. 
Let $v \in \pesos$. We claim that  the sequence $\cV_v  = (v_i \, , \, V_i)_{i\in \ J}$ 
is a FS refinement of $\cW_v\,$. 

Indeed, consider $T_{\cV_v} = T_{\cW_v}\big|_{\cK_\cV}\,
\in L(\cK_\cV \, , \, \H)$. By Lemma \ref{proy}, $\dim \cK_\cV\orto = 
\exe{\cW , \cV}<\infty$. As in the proof of Lemma \ref{exari}, this implies that 
$R(T_{\cV_v} ) = R(T_{\cW_v }\, P_{\cK_\cV} ) \subcer \H$. On the other hand, 
$\sspan\{ \ds\cup _{i \in J} V_i\} \inc R(T_{\cV_v} )$. But $\sspan\{ \ds\cup _{i \in J} V_i\}$ 
is dense in $\H$, because $T_{\cV_w}$ is surjective (recall that $\cV_w$ is a FS). 
This shows that also $T_{\cV_v}$ is surjective, i.e. $\cV_v  $  is a FS as claimed. 
In other words, we have that $\cV$ is a RBS, and $v_J= \{v_i\}_{i\in \ J} \in \cP(\cV)$. 
By Proposition \ref{los pesos}, 
$v_J \in \linmJ^*$. As before, this implies that $v \in \bsi$. 
Using Proposition \ref{los pesos} again, we 
conclude that $\exe {\cW_v } = \exe {\cW_w } $.  \QED

\begin{teo} \label{raros2}
Let $\sfram$ be a FS for $\hil$. 
Then 
$$
\exe {\cW_v } = \exe {\cW_w } \peso{for every other} v \in \pesos  \ .
$$
\end{teo}
\proof If  $\exe {\cW_w } < \infty$, apply  Corollary \ref{raros}. 
If  $\exe {\cW_w } = \infty$ and $v \in \pesos$, 
then  also $\exe {\cW_v }= \infty$, since otherwise we could apply 
Corollary \ref{raros} to $\cW_v\,$.  \QED

\section{Examples}
Observe that, if $\{ E_i\}_{i \in I}$ is an OBS of $\cK$ and  $T\in L(\cK , \hil)$
is a surjective operator such that $T(E_i ) \subcer \H$ for every $i \in I$, 
then  $\cW = \{ TE_i\}_{i \in I}\, $ is a generating sequence for $\H$. 
Nevertheless, our first example shows that, in general, such a sequence $\cW $
may have $\pesos = \vacio$, i.e. $\FS$  
fails to be a FS for $\hil$, for any sequence $w\in \linm$ of weights.

\begin{exa}\label{gama no frame}
Take $\cB = \{e_n\}_{n\in \NN}$ an \onb of $\hil$. For every 
for $k\in \NN$, consider the space $E_k=\sspan\{e_{2k-1},e_{2k}\}\,$. 
Observe that   $E_k$ is an OBS for $\hil$.
Consider  the (densely defined) operator $T: \hil \to \hil$  given by 
\[
Te_{n}= \begin{cases} 
2^{-k} \, e_1 
& \mbox{if $n = 2k-1 $}\\ & \\ 
e_{k+1}                 &\mbox{if $n = 2k$} 
\end{cases} \quad .
\]
Then, $T$ can be extended to a bounded surjective operator $T$, since the sequence
$\{Te_k\}_{k\in \NN}$ is easily seen to be a tight frame for $\hil$. 
We shall see  that the sequence of closed subspaces 
$$
\cW = \{W_k\}_{k\in \N} \peso{given by} W_k=T(E_k)=\sspan\{e_1,e_{k+1}\}  \ \ , \quad k \in \NN 
$$ 
satisfies that $\pesos = \vacio$. 
Indeed, suppose that $w \in \pesos$. 
Then by Eq. \eqref{frame de sub} applied to 
$\displaystyle  f=e_1\in \displaystyle \bigcap_{k\in \NN}W_k\, $, we  would have 
that $w\in \ell\,^2(\NN )$. But this contradicts the existence of a 
lower frame bound $A_{\cW_w}  $ for  $\sframN \,$, because for every $k \in \NN$, 
 $$
A_{\cW_w}   = A_{\cW_w}  \|e_{k+1}\|^2  \le 
\sum_{j\in\NN} w_j^2\|P_{W_j}e_{k+1}\|^2=w_k^2 \xrightarrow[k\ra \infty]{} 0 \ .
 $$
Observe that, by definition,  
$\frac{\ds \gamma(TP_{E_k})}{\ds\|TP_{E_k}\|} = \frac{\ds 2^{-k}}{\ds 1}
 \convk 0$.  
\EOE
\end{exa}

\medskip
\noi
The operator $T$ and the OBS $\cE = \{E_n\}_{k\in \NN}$ of the last Example do 
not satisfy  Eq. \eqref{cota gama} in Theorem \ref{cotas}. 
Still, Eq.  \eqref{cota gama}  is not a necessary  condition in order 
to assure  that $\pesos \neq \vacio $, if $\cW = T\cE$. 
Next example shows a FS wich is the image of an 
OBS under an epimorphism which does not satisfy Eq. \eqref{cota gama}.

\begin{exa}\label{no dos}
Let $\bon$ be an orthonormal   basis for $\hil$ and consider the   frame (of vectors) 
\[
\F= \fram \peso{given by } 
f_{n}= \begin{cases}  e_k 
& \mbox{if $n = 2k-1 $}\\  & \\ 
\frac{\ds e_{k+1}}{\sqrt{k+1}}   &\mbox{if $n = 2k$} 
\end{cases} \quad .
\]
Let $T=T_\F \in L(\ell^2(\NN ),  \hil)$ be its synthesis operator (wich is surjective). 
If  $\{b_n\}_{n\in \NN}\,$ is 
the canonical  basis of $\ell^2(\NN )$, then $Tb_n=f_n\, $. For each $k\in \NN$ we set $E_k=\sspan\{b_{2k-1},b_{2k}\}$. Then, by construction, $\{E_k\}_{k\in \NN}$ is an 
OBS of $\ell^2(\NN )$.  
Take the sequences 
$$
w =  e \in \linmN \peso{and}   W_k = TE_k=\sspan\{e_k,e_{k+1}\} \quad , \quad k \in \NN \  .
$$ 
By Theorem \ref{CK},  $\sframN $  is a  FS for $\hil$.
Nevertheless, $T$ does not satisfy Eq. \eqref{cota gama}, since 
$\gamma(TP_{E_k})=\frac{\ds 1}{\sqrt{k+1}}\,$, 
while $\|TP_{E_k}\| = 1 $, for every $k \in \NN$.
\EOE
\end{exa}

\medskip
\noi
The key argument in Example \ref{gama no frame} 
was  that  $\ds \bigcap_{i\in I}W_i \neq \trivial$. This fact is 
sufficient for the emptiness of  $\pesos$ if $\sspan \{W_i : 1\le i \le n\} \neq \hil$ 
for every $n \in \NN$. Nevertheless, next example shows  
a minimal and generating sequence $\cW$  of  finite dimensional subspaces 
such that $\pesos = \vacio$.

\begin{exa} \label{sin int}
Fix an \onb  $\cB = \{e_i\}_{i\in \NN}$  for 
$\hil$. Consider the unit vector $\ds g =\sum_{k=1}^\infty \frac{e_{2k}}{2^{k/2}} \in \hil$. 
For every $n\in \NN$, denote by $P_n \in \op$ the orthogonal projection onto  
$\hil_n = \sspan \{e_1,e_2,\ldots,e_n\}$. Consider the generating sequence 
$\cW= \{W_k\}_{k\in \NN}\,$ given by  
$$
W_k=\sspan\{P_{2k}\ g\, , \ e_{2k-1}\} = \sspan \left\{\sum_{j=1}^k \frac{e_{2j}}{2^{j/2}} \ , \  e_{2k-1}\right\}
   \quad , \quad k \in \NN \ .
$$
Straightfordward computations show that $\cW$ is a minimal sequence. 
The problem is that 
$\angf{W_i\,}{W_j} \xrightarrow[i,j \rightarrow\infty]{} 0$
 exponentially, and for this reason $\pesos = \vacio$.  Indeed, suppose that $w \in \pesos$, 
 and that  $\cW_w = (w, \cW)$ is a FS.  Then 
\begin{align}  
 B_{\cW_w}= B_{\cW_w} \|g\|^2 
  \ge \sum_{k\in \NN}\ w_k ^2 \ \|P_{W_k}\, g\|^2&=\sum_{k\in \NN}\ w_k ^2 \, \|P_{2k}\, g\|^2=
\sum_{k\in \NN}\ w_k ^2 \, (1-2^{-k}) 
\ ,
\end{align}
which implies that  $w_k \xrightarrow[k\ra \infty]{} 0$. On the other hand, for every $k\in \NN$,
\begin{equation}
A_{\cW_w} = A_{\cW_w}\|e_{2k-1}\|^2 \le \sum_{i\in \NN}\ w_i ^2 \, \|P_{W_i}\ e_{2k-1}\|^2=w_k ^2 
 \ \ \implies \ \ A_{\cW_w}=0  \ , 
\end{equation}
a contradiction. So $\pesos = \vacio$. \EOE
\end{exa}

\medskip
\noi
It is well known that $\{f_j\}_{j\in \NN}$ is a Parseval frame in $\hil$ if and only if 
there exists a Hilbert $\cK$ containing $\hil$ such that $f_j=P_{\hil}b_j$ for every 
$j\in \NN$, where $\{b_j\}_{j\in \NN}$ is an \onb for $\cK$. One may think that 
a similar result is true for tight frames of subspaces, where we replace orthonormal 
basis by OBS. In section 4 we proved one implication (a Parseval
FS is an orthogonal projection of an OBS) but the 
converse it is not true: 

\begin{exa}\label{proyeccion no frame}
Let $\bon$  be an orhonormal basis for $\hil$. Consider the unit vector 
$$
g=\ds \sum_{k\in \NN} \frac{e_{2k-1}}{2^{k/2}} \ , 
\peso{and take}  \cM=\overline{\sspan}\, \{g \} \cup \{e_{2k} : k\in \NN\} \ .
$$ 
On the other hand, take the sequence $\cE = \{E_k\}_{k\in \NN}$ given by 
$E_k=\sspan\{ e_{2k-1},e_{2k}\}$ ($k \in \NN$). Then  
$\cE$ is an OBS for $\hil$.
Take the sequence 
$$
\sseck  \peso{given by}  W_k=P_\cM  \, E_k=\sspan\{g\, , \, e_{2k}\} \ ,  \peso{for every $k \in \NN$ .} 
$$  
Then $\pesos = \vacio$   
by same reason as in  Example \ref{gama no frame}, 
because $\ds g \in \bigcap_{k \in \NN} W_k \neq \trivial $  . 
\EOE
\end{exa}

\begin{exa} \label{frame no gama}\rm 
Let $\cE = \{e_n\}_{n\in \NN}$ be an \onb of $\hil$. 
Consider the sequence $\sseck \,$ given by 
$$
W_1 = \gen { e_k : k \ge 2 } = \{e_1\}\orto \peso{and}
W_k  = \sspan\{e_1,e_{k }\}  \ \ , \peso{for} k \ge 2  \  .
$$ 
Observe that 
$\pesos  = \ell_+^{\,2}(\NN) $. Indeed,  one inclusion is clear,  and  
$$ 
w \in \pesos \implies 
\sum _{k=2}^\infty  w_k^2 \ = \ \sum_{k=2}^\infty   w_k^2\ \|P_{W_k}e_1\|^2
\le B_{\cW_w} \ \implies \ w \in \ell_+^{\ 2}(\NN)  \ . 
$$ 
Now we shall see that   $\cW_w $ can not be a {\bf tight} FS for any 
$w \in \pesos$. Indeed, 
if $\cW_w$ 
where a $A$-tight frame, then 
for every $k\ge 2 $,
$$\ds 
A = A\| e_k\|^2 = \sum_{i\in \NN} w_i^2 \ \|P_{W_i}e_k\|^2=w_1^2+w_{k}^2 
 \ \ \implies \ \  w_{k}^2  = A -w_{1}^2   \ , $$ 
which contradicts the fact that $w \in \ell_+^{\,2}(\NN)$.
Our next step is to show that the frame operator $S_{\cW_w}\in L(\hil)$ is diagonal 
with respect to $\cE$, for every $w \in  \pesos$. Indeed, 
$$
\tfs^* e_1=\{w_kP_{W_k}e_1\}_{k\in \NN}={0}\oplus \{w_ke_1\}_{k\geq 2} \implies 
S_{\cW_w}e_1=\tfs \tfs^* e_1=\left(\sum_{k=2}^\infty  w_k^2\right)e_1 \ .
$$
On the other hand, if $E_k $ is the copy of each $W_k $ in $\cK_{\cW}\,$, then for every  
$k \in \NN$ and $j\ge 2$, 
$$P_{E_k} \left(\tfs^*e_j\right)=\begin{cases} 
w_1 \ e_j 
& \mbox{if $k=1 $}\\   
w_{j}\ e_j   &\mbox{if $k=j$}\\
\ \ 0 &\mbox{if $k\neq 1, j$} 
\end{cases} \implies 
S_{\cW_w} e_j=\tfs \tfs^* e_j=(w_1^2+w_{j}^2)e_j \ .
$$
In particular,  $S_{\cW_w}\mrai$ is also diagonal. This implies that 
$S_{\cW_w}\mrai \cW  = \cW $,  which we have seen that can not be  tight for any sequence of weights.

Another property of this example is the following: $\FS$ is a FS for $\H$, but 
the sequence $(w_k \, , \, W_k)_{k>1}$ is not a frame sequence of subspaces 
(i.e. a FS for $\gen{W_k : k>1}$). This can be proved by the same argument as 
in Example \ref{gama no frame}, using that $\cap _{k>1} W_k \neq \trivial$. 
\EOE
\end{exa}

\begin{exa} \label{otros2}\rm 
Let $\cB_4 = \{e_n\}_{n\le 4}$ be an \onb of $\C^4$. 
Consider the sequence 
$$
W_1 = \sspan \{ e_1 , e_2\} \ , \quad  W_2  = \sspan \{ e_1 , e_3\} \peso{and}
W_3  = \sspan\{e_{4 }\}  \ .
$$ 
We shall see that, for every invertible 
$G \in \cM_{4}(\C )$, 
and every $w \in \R_+^3\,$, the sequence
$G\FS = (w_k , GW_k)_{k \in \IN{3}} $ fails to be a Parseval FS. 
Take \onb of each $GW_i$ 
$$
GW_1 = \sspan \{ g_1, g_2 \} \ , \quad GW_2 = \sspan \{ g_1, g_3 \} \peso{and} 
GW_3 = \sspan \{ g_4 \} \ ,
$$
where $\ds g_1 = \frac{Ge_1}{\|Ge_1\|}\ $, and similarly for $g_4\,$. 
If $G\FS$ were a Parseval FS, then 
the frame  
$$\cE = \{T_ {G\FS} g_k\}_{k\in \IN{5}}\,
= \{ w_1\  g_1 \, ,\,  w_1\  g_2\,  ,\,  w_2 \ g_1\,  ,\,  
w_2 \ g_3\,  ,\,  w_3 \ g_3\} \ , 
$$
would be also Parseval. Consider the matrix $T \in \cM_{4,5}(\C )$ with the vectors of 
$\cE$ as columns. After a unitary change of coordinates, 
$T$ has the form 
$$
T = \bm{ccc} w_1 & w_2 & \vec{v}  \\0&0& V \em \barr{l} \C \\ \C^3 \earr 
\peso{with $\vec{v} = (0,0,a)\in \C^3$ and $V \in \cM_3(\C)$  .}
$$
Since $TT^* =I_4\,$, it is easy to see that $V\in \cU(3)$. But this is 
impossible because the first two columns of $V$ have norms 
$\|w_1\, g_2\| = w_1$ and $\|w_2\, g_3\|= w_2\,$, while $1 = w_1^2 +w_2^2 +|a|^2\,$. \EOE
\end{exa}

\noi
{\bf Acknowledgments}: We wish to thank Professors G. Kutyniok and P. Gavruta for bringing 
to our attention their recently works on  fusion frames 

\fontsize {9}{10}\selectfont

\vglue0.5truecm
\noi {\bf Mariano A. Ruiz and Demetrio Stojanoff
}

\noi Depto. de Matem\'atica, FCE-UNLP,  La Plata, Argentina
and IAM-CONICET  

\noi e-mail: maruiz@mate.unlp.edu.ar and demetrio@mate.unlp.edu.ar

\end{document}